\documentclass[10pt]{article}
\usepackage[latin1]{inputenc}
\usepackage[T1]{fontenc}
\usepackage{amssymb,amstext}
\usepackage[english]{babel}

\usepackage{amsmath,amsthm}
\usepackage{graphicx,graphics,here}

\renewcommand{\(}{\left(}
\renewcommand{\)}{\right)}
\renewcommand{\[}{\left[}
\renewcommand{\]}{\right]}

\newcommand{\F}{{\mathcal F}}
\renewcommand{\L}{{\mathcal L}}

\newcommand{\R}{\mathbb R}
\newcommand{\C}{\mathbb C}
\newcommand{\N}{\mathbb N}

\renewcommand{\P}{\mathrm P}
\newcommand{\E}{\mathrm E}

\newcommand{\partI}[2]{\frac{\partial#1}{\partial#2}}
\newcommand{\partII}[2]{\frac{\partial^2#1}{\partial#2^2}}

\newtheorem{theorem}{Theorem}[section]
\newtheorem{lemma}[theorem]{Lemma}
\newtheorem{proposition}[theorem]{Proposition}
\newtheorem{corollary}[theorem]{Corollary}

\theoremstyle{definition}
\newtheorem{definition}[theorem]{Definition}

\theoremstyle{remark}
\newtheorem{remark}[theorem]{Remark}

\numberwithin{equation}{section}

\title{Towards a characterization of Markov processes enjoying the
time-inversion property}
\author{Stephan Lawi\footnote{Laboratoire de Probabilit\'es et Mod\`eles Al\'eatoires,
CNRS-UMR 7599, Universit\'e Paris VI \& Universit\'e Paris VII, 4
Place Jussieu, 75252 Paris Cedex 05, France. \textit{Email}:
lawi@proba.jussieu.fr.}}

\date{\today}

\begin{document}

\maketitle

\begin{abstract}
We give a necessary and sufficient condition for a homogeneous
Markov process taking values in $\R^n$ to enjoy the time-inversion
property of degree $\alpha$. The condition sets the shape for the
semigroup densities of the process and allows to further extend the
class of known processes satisfying the time-inversion property. As
an application we recover the result of Watanabe in \cite{Wa1975}
for continuous and conservative Markov processes on $\R_+$. As new
examples we generalize Dunkl processes and construct a matrix-valued
process with jumps related to the Wishart process by a skew-product
representation.
\end{abstract}
\vspace{10pt}

\\



\section{Introduction}
Let $\{(X_t,t\ge0);\ (\P_x)_{x\in\R^n}\}$ be a homogeneous Markov
process with semigroup densities (assumed to exist and to be twice
differentiable in the space and time variables):
\begin{equation}
P_t(x,dy) = p_t(x,y) dy.
\end{equation}
For all $x\in\R^n$ and some $\alpha>0$, the process $\{(t^\alpha
X_\frac{1}{t},t>0);\ \P_x\}$ is Markov and in general inhomogeneous.

\begin{definition}
The process $\{(X_t,t\ge0);\ \P_x\}$  is said to enjoy the {\it
time-inversion property of degree $\alpha$} if the Markov process
$\{(t^\alpha X_\frac{1}{t},t>0);\ \P_x\}$ is homogeneous.
\end{definition}

The problems associated to time-inversion of Markov processes are
closely related to the so-called dual processes in probabilistic
potential theory (see \cite{Dy1975} and \cite{Ku1977}). A pair of
dual processes is a pair of Markov processes whose resolvents are
conjugate relative to some measure. The trajectories of these
processes are connected by means of time-reversal. However,
time-reversal in general violates homogeneity and changes the path
properties of the process. In this paper, we are mainly concerned
with the first issue, and in particular, whether one can
characterize the semigroup densities of a Markov process such that
homogeneity is preserved under time-inversion.

Celebrated examples of Markov processes, known to enjoy the
time-inversion property for $\alpha=1$, are Brownian motions with
drift in $\R^n$ and Bessel processes with drift (see \cite{PY1981}
and \cite{Wa1975}). Gallardo and Yor \cite{GY2005b} recently worked
out a sufficient condition on the semigroup densities for a Markov
process to enjoy the time-inversion property. Their argument
extended the class of processes to processes with jumps such as the
Dunkl process \cite{RV1998} and matrix-valued processes such as the
Wishart process \cite{Br1991}. The aim of the paper is to find a
necessary and sufficient condition and to provide some new examples.

Section \ref{sec:def} contains the main theorem of the paper, which
is proved in section \ref{sec:proof} using straightforward
analytical arguments. Section \ref{sec:applications} considers an
application of the theorem to Markov processes on $\R_+$. The result
is shown to be strong enough to entirely characterize the class of
diffusion processes on $\R_+$ that enjoy the time-inversion property
and thus provide a different proof than that of Watanabe in
\cite{Wa1975}. Section \ref{sec:examples} gives new examples of
processes that enjoy the time-inversion property. We review how the
generalized Dunkl process fits the requirements of the theorem and
then introduce a matrix-valued process with jumps. The relation of
the latter process to the Wishart process mimics the relation
between the one-dimensional Dunkl and Bessel processes.

\section{Markov processes which enjoy time-inversion}\label{sec:def}
Fix $x\in\R^n$. Recall that under $\P_x$, the process $(t^\alpha
X_\frac{1}{t},t>0)$ is Markov, inhomogeneous, with transitional
probability densities $q_{s,t}^{(x)}(z,y),\ (s<t;z,y\in\R^n)$, which
satisfy the following relation:
\begin{equation}
E_x\[f(t^\alpha X_\frac{1}{t})\ \big|\ s^\alpha X_\frac{1}{s}=z\] =
\int dy\ f(y)\ q_{s,t}^{(x)}(z,y)
\end{equation}
where
\begin{equation}\label{eq:qp}
q_{s,t}^{(x)}(a,b) = t^{-n\alpha}\
\frac{p_\frac{1}{t}\(x,\frac{b}{t^\alpha }\)\
p_{\frac{1}{s}-\frac{1}{t}}\(\frac{b}{t^\alpha },\frac{a}{s^\alpha
}\)} {p_\frac{1}{s}\(x,\frac{a}{s^\alpha }\)}.
\end{equation}

The process $(t^\alpha X_\frac{1}{t},t>0)$ is not uniquely defined
from the knowledge of the semigroup densities $p_t(x,y)$. Actually,
there exists at least one transformation that leaves the semigroup
densities $q_{s,t}^{(x)}(a,b)$ unchanged: Doob's $h$-transform.

\begin{definition}
Doob's $h$-transform is the transformation
\begin{equation}
T_h: \P_x|_{\F_t} \mapsto \frac{h(X_t)}{h(x)}e^{-\nu t}\
\P_x|_{\F_t}.
\end{equation}
for some function $h$ and some constant $\nu>0$.
\end{definition}

This remark leads to our first assertion.

\begin{proposition}\label{prop:h}
Two processes related by $h$-transforms yield the same process by
time-inversion.
\end{proposition}
It is hence legitimate to research for a criterium to classify all
processes that enjoy the time-inversion property up to
$h$-transforms. The following theorem gives a concise statement of
our main result:
\begin{theorem}\label{thm:main}
The Markov process $(t^\alpha X_\frac{1}{t}, t>0)$ is homogeneous if
and only if the semigroup densities of $(X_t,t\ge0)$, assumed twice
differentiable, are of the form:
\begin{equation}\label{eq:density}
p_{t}(x,y) = t^{-\frac{n\alpha}{2}}\
\Phi\(\frac{x}{t^\frac{\alpha}{2}}, \frac{y}{t^\frac{\alpha}{2}}\)\
\theta\(\frac{y}{t^\frac{\alpha}{2}}\)\
\exp\bigg\{\rho\(\frac{x}{t^\frac{\alpha}{2}}\)
+\rho\(\frac{y}{t^\frac{\alpha}{2}}\)\bigg\},
\end{equation}
or if they are in $h$-transform relationship with it. The functions
$\Phi, \theta, \rho$ have the following properties for $\lambda>0$:
\begin{enumerate}
\item[1.] $\Phi(\lambda x,y) = \Phi(x,\lambda y)$;

\item[2.] $\rho(\lambda x) = \lambda^\frac{2}{\alpha}\ \rho(x)$;

\item[3.] $\displaystyle\theta(\lambda y) =
\lambda^{\beta}\ \theta(y)$ for some $\beta\in\R$;

\end{enumerate} Moreover, if the symmetry condition $\Phi(x,y) =
\Phi(y,x)$ is satisfied, then the semigroup densities are related as
follows:
\begin{equation}\label{eq:qhp}
q_{t}^{(x)}(a,b) = \frac{\Phi\(x,b\)} {\Phi\(x,a\)}\
\exp\Big\{t\rho\(x\)\Big\}\ p_t(a,b).
\end{equation}
\end{theorem}
Such a decomposition of the semigroup densities shows furthermore
that $h$-transforms are the only transformations that leave
semigroup densities of the time-inverted process $(t^\alpha
X_\frac{1}{t},t>0)$ unchanged.

\begin{definition}
A Markov process $(X_t,t\ge0)$ is called semi-stable with index
$\gamma$ in the sense of Lamperti \cite{La1972} if:
\begin{equation}
\Big\{(X_{ct}, t \ge0);\P_x\Big\} \stackrel{(d)}{=} \Big\{(c^\gamma
X_{t}, t \ge0);\P_{x/c^\gamma}\Big\}.
\end{equation}
\end{definition}
\noindent As a consequence, the semigroup densities of a semi-stable
process with index $\gamma$ have the following property:
\begin{equation}
p_t(x,y) = t^{-n\gamma}\ p_1\(\frac{x}{t^\gamma}
,\frac{y}{t^\gamma}\).
\end{equation}
Remark that the expression for the semigroup densities in Theorem
\ref{thm:main} satisfies this property for $\gamma =\alpha/2$. This
remark yields the following corollary:

\begin{corollary}\label{cor:main}
A Markov process that enjoys the time-inversion property of degree
$\alpha$ is a semi-stable process of index $\alpha/2$, or is in
$h$-transform relationship with it. The converse is not true.
\end{corollary}

\begin{remark}
Throughout the paper, we assume the semigroup densities $p_t(x,y)$
to be positive over the domain and regular enough to be at least
twice differentiable in the space and time variables. This
assumption is more a technical requirement for the proof than a
necessity for the characterization of Markov processes enjoying the
time-inversion property.
\end{remark}

\section{Proof of the Theorem}\label{sec:proof}

\subsection{Sufficiency}
If $p_t(x,y)$ satisfies the condition (\ref{eq:density}), then from
formula (\ref{eq:qp}) the semigroup densities $q_{s,t}^{(x)}(a,b)$
can be written as
$$
q_{s,t}^{(x)}(a,b) = (t-s)^{-\frac{n\alpha}{2}}\ \prod_{i=1}^3
R_{s,t}^{(i)}
$$
where
\begin{eqnarray*}
R_{s,t}^{(1)} &=& \frac{\Phi\(x\ t^\frac{\alpha}{2},
\frac{b}{t^\alpha}\ t^\frac{\alpha}{2}\)} {\Phi\(x\
s^\frac{\alpha}{2},\frac{a}{s^\alpha}\ s^\frac{\alpha}{2}\)}\
\Phi\(\frac{b}{t^\alpha}\(\frac{st}{t-s}\)^\frac{\alpha}{2},
\frac{a}{s^\alpha}\(\frac{st}{t-s}\)^\frac{\alpha}{2}\)
\\
R_{s,t}^{(2)} &=& \frac{\theta\(\frac{b}{t^\alpha}\
t^\frac{\alpha}{2}\)} {\theta\(\frac{a}{s^\alpha}\
s^\frac{\alpha}{2}\)}\
\theta\(\frac{a}{s^\alpha}\(\frac{st}{t-s}\)^\frac{\alpha}{2}\) \\
R_{s,t}^{(3)} &=& \exp\left\{\rho\(x\ t^\frac{\alpha}{2}\) +
\rho\(\frac{b}{t^\alpha}\ t^\frac{\alpha}{2}\)\ +
\rho\(\frac{b}{t^\alpha}\(\frac{st}{t-s}\)^\frac{\alpha}{2}\)\right.\\
&&  \quad\quad +\
\rho\(\frac{a}{s^\alpha}\(\frac{st}{t-s}\)^\frac{\alpha}{2}\)-
\rho\(x\ s^\frac{\alpha}{2}\) - \rho\(\frac{a}{s^\alpha}\
s^\frac{\alpha}{2}\)\Bigg\}
\end{eqnarray*}
Using the properties of $\Phi,\theta,\rho$ described in Theorem
\ref{thm:main}, we obtain:
\begin{eqnarray*}
R_{s,t}^{(1)} &=& \frac{\Phi\(x,b\)} {\Phi\(x,a\)}\
\Phi\(\frac{b}{(t-s)^\frac{\alpha}{2}},
\frac{a}{(t-s)^\frac{\alpha}{2}}\)
\\
R_{s,t}^{(2)} &=& (t-s)^{-\frac\alpha2\beta}\ \theta(b)\
\\
R_{s,t}^{(3)} &=& \exp\left\{(t-s)\rho(x) +
\rho\(\frac{a}{(t-s)^\frac{\alpha}{2}}\) +
\rho\(\frac{b}{(t-s)^\frac{\alpha}{2}}\)\right\}.
\end{eqnarray*}
Hence there is no separate dependence on $s$ and $t$, but only on
the difference $t-s$, which allows to conclude that
$$
q_{s,t}^{(x)}(a,b) = q_{t-s}^{(x)}(a,b)
$$
and proves homogeneity for the process $(t^\alpha
X_\frac{1}{t},t>0)$. If in addition $\Phi(x,y) = \Phi(y,x)$, then
\begin{eqnarray*}
q_{s,t}^{(x)}(a,b) &=& \frac{\Phi\(x,b\)}{\Phi\(x,a\)}\
\frac{e^{(t-s)\rho(x)}}{(t-s)^{\frac{n\alpha}{2}}}\
\Phi\(\frac{a}{(t-s)^\frac{\alpha}{2}},
\frac{b}{(t-s)^\frac{\alpha}{2}}\)\
\theta\(\frac{b}{(t-s)^{\frac\alpha2}}\)\\
&& \exp\left\{\rho\(\frac{a}{(t-s)^\frac{\alpha}{2}}\) +
\rho\(\frac{b}{(t-s)^\frac{\alpha}{2}}\)\right\},
\end{eqnarray*}
in which we recognize (\ref{eq:qhp}).

\subsection{Necessity}
For simplicity of notation, we prove the necessity of condition
(\ref{eq:density}) in the case $\alpha=1$. The extension to
$\alpha>0$ is immediate by the change of variables $x\mapsto
\|x\|^\alpha\frac{x}{\|x\|}$. Recall first the following definition
of homogeneous functions:
\begin{definition} A function
$f:\R^n\to\R$ that satisfies for $x=(x_1,\ldots,x_n)$,
$\forall\lambda\in\R_+$,
$$
f(\lambda x) = \lambda^\beta f(x)
$$
is called {\it homogeneous of degree $\beta\in\R$}. If $f\in
C^1(\R^n)$, Euler's Homogeneous Function Theorem gives a necessary
and sufficient condition for the function $f(x)$ to be homogeneous:
$$
\sum_{i=1}^n x_i \partI{}{x_i} f(x) = \beta f(x).
$$
\end{definition}

Consider the function of $2n+1$ variables: $l(x,y,t) =
\ln(p_t(x,y)),\ x,y\in\R^n$. From (\ref{eq:qp}), $l$ must satisfy
for $s=t-h$:
\begin{equation}\label{eq:ldt}
\partI{}{t}\[n\ln\frac{1}{t} + l\(x,\frac{b}{t},\frac{1}{t}\)
+ l\(\frac{b}{t},\frac{a}{t-h},\frac{1}{t-h}-\frac{1}{t}\) -
l\(x,\frac{a}{t-h},\frac{1}{t-h}\)\] =0.
\end{equation}

\subsubsection{The kernel $\Phi(x,y)$}
Taking derivatives with respect to $b_i$ and $a_j$ for some $0\le
i,j\le n$ yields:
$$
\partI{}{b_i}\partI{}{a_j}\partI{}{t}\
l\(\frac{b}{t},\frac{a}{t-h},\frac{1}{t-h}-\frac{1}{t}\)=0.
$$
If we set $\phi\(x,y,t\) = \partI{}{x_i}\partI{}{y_j} l\(x,y,t\)$,
the latter becomes
$$
-\frac{1}{t^2(t-h)}\(\phi + \frac{b}{t}\cdot\nabla_1\phi\) -
\frac{1}{t(t-h)^2}\(\phi + \frac{a}{t-h}\cdot\nabla_2\phi\)$$
$$
+ \frac{1}{t(t-h)}\(\frac{1}{t^2}-\frac{1}{(t-h)^2}\)\partial\phi=0,
$$
with the notation $\nabla_1 =
(\partI{}{x_1},\ldots,\partI{}{x_n})^T,\ \nabla_2 =
(\partI{}{y_1},\ldots,\partI{}{y_n})^T$, $\partial$ the derivative
with respect to the time variable and $\phi =
\phi\(\frac{b}{t},\frac{a}{t-h},\frac{1}{t-h}-\frac{1}{t}\)$. For
clarity, we change variables to
$$
z_1 = \frac{b}{t}\quad,\quad z_2=\frac{a}{t-h}\quad,\quad
t_1=-\frac{1}{t}\quad,\quad t_2=\frac{1}{t-h}.
$$
Then $\phi = \phi(z_1,z_2,t_1+t_2)$ and
$$
t_1\bigg(\phi + z_1\cdot\nabla_1\phi\ + t_1\ \partial\phi\bigg) =
t_2\bigg(\phi + z_2\cdot\nabla_2\phi\ + t_2\ \partial\phi\bigg),
$$
or equivalently
$$
t_1\bigg(\phi + z_1\cdot\nabla_1\phi\ + (t_1+t_2)\
\partial\phi\bigg) = t_2\bigg(\phi + z_2\cdot\nabla_2\phi\ + (t_1+t_2)\
\partial\phi\bigg).
$$
We change variables once more: $u=\frac{t_1}{t_2}, v=t_1+t_2$ to
get:
\begin{equation}\label{eq:phi}
\phi + z_1\cdot\nabla_1\phi\ + v\
\partial\phi = u\big(\phi + z_2\cdot\nabla_2\phi\ + v\
\partial\phi\big).
\end{equation}
This gives the following proposition.

\begin{proposition}
Equation (\ref{eq:phi}) is satisfied if and only if
\begin{equation}\label{eq:phihomog}
\phi(\lambda x,\mu y,\lambda\mu t) = \frac{1}{\lambda\mu}\
\phi(x,y,t).
\end{equation}
\end{proposition}

\begin{proof}
As the LHS of equation (\ref{eq:phi}) is independent of $u$, one can
readily take the equivalent condition:
\begin{eqnarray}\label{eq:phi01}
\phi + z_1\cdot\nabla_1\phi\ + v\
\partial\phi &=& 0,\\
\phi + z_2\cdot\nabla_2\phi\ + v\
\partial\phi &=& 0.\label{eq:phi02}
\end{eqnarray}
Let $g(\lambda) = \phi(\lambda z_1,z_2,\lambda v)$ and $h(\mu)
=\phi(z_1,\mu z_2,\mu v)$. Equation (\ref{eq:phi01}) implies $
g(\lambda) +\lambda g'(\lambda) =0,$ which solves to $g(\lambda) =
\lambda^{-1}g(1)$ and hence
$$
\phi(\lambda z_1,z_2,\lambda v) = \frac{1}{\lambda}\phi(z_1,z_2,v).
$$
Equation (\ref{eq:phi02}) implies $ h(\mu) +\mu h'(\mu) =0,$ which
solves to $h(\mu) = \mu^{-1}h(1)$ and hence
$$
\phi(z_1,\mu z_2,\mu v) = \frac{1}{\mu}\phi(z_1,z_2,v).
$$
Combining the latter two equations yields (\ref{eq:phihomog}).

Conversely, if $g$ is homogeneous of degree $-1$ and $h$ is
homogeneous of degree $-1$, we get:
$$
g(\lambda) + \lambda g'(\lambda) = h(\mu) + \mu h'(\mu) = 0,
$$
which is equivalent to (\ref{eq:phi01}) and (\ref{eq:phi02}) and
concludes the proof.
\end{proof}

The scaling property (\ref{eq:phihomog}) implies moreover the
equivalent formulation
\begin{equation}\label{eq:phi1t}
\phi(x,y,t) = \frac{1}{t}\phi_1\(\frac{x}{\sqrt t}, \frac{y}{\sqrt
t}\),
\end{equation}
where $\phi_1(x,y) = \phi(x,y,1)$ must satisfy
\begin{equation}\label{eq:phi1homog}
\phi_1(\lambda x,y) = \phi_1(x, \lambda y).
\end{equation}
Under the change of variables $\displaystyle\bar x = \frac{x}{\sqrt
t}$, $\displaystyle\bar y = \frac{y}{\sqrt t}$, the kernel $\bar
l(\bar x,\bar y,t) = l(x,y,t)$ satisfies
$$
\phi_1(\bar x,\bar y) = \partI{}{\bar x_i}\partI{}{\bar y_j}\bar
l(\bar x,\bar y,t).
$$
Note that the previous results remain valid for all second
derivatives of $\bar l(\bar x,\bar y,t)$ with respect to $\bar x_i$
and $\bar y_j$, $\forall i,j\in\{1,\ldots,n\}$. By integration over
$\bar x_i$ and $\bar y_j$, one can thus already make an assumption
on the general shape of the kernel $\bar l(\bar x,\bar y,t)$:
$$
\bar l(\bar x,\bar y,t) = k(\bar x,\bar y) + \bar\varphi_1(\bar x,t)
+ \bar\varphi_2(\bar y,t),
$$
where $k:\R^{2n}\to\R$ does not depend explicitly on time and
$\bar\varphi_1,\bar\varphi_2:\R^{n+1}\to\R$. An equivalent
representation for $l(x,y,t)$ gives
\begin{equation}\label{eq:generalforml}
l(x,y,t) = k\(\frac{x}{\sqrt t},\frac{y}{\sqrt t}\) + \varphi_1(
x,t) + \varphi_2(y,t),
\end{equation}
with $\bar\varphi_1(\bar x,t) = \varphi_1(x,t)$ and
$\bar\varphi_2(\bar y,t) = \varphi_2(y,t)$.

The scaling property (\ref{eq:phi1homog}) of $\phi_1(x,y)$ now
translates to $k(x,y)$ as follows:
\begin{equation}\label{eq:kscaling0}
\partI{}{x_i}\partI{}{y_j}k(\lambda x,y) =
\partI{}{x_i}\partI{}{y_j}k(x,\lambda y),\quad
\forall i,j \in \{1,\ldots,n\}.
\end{equation}
By integration, $k(x,y)$ must satisfy the following condition,
\begin{equation}\label{eq:kscaling1}
k(\lambda x, y) = k(x,\lambda y) + \zeta_1(x,\lambda)+
\zeta_2(y,\lambda),
\end{equation}
along with the properties $\zeta_1(x,1) + \zeta_2(y,1)=0$ and
$$
\zeta_1(x,\alpha\lambda) + \zeta_2(y,\alpha\lambda)  =
\zeta_1(x,\alpha) + \zeta_1(\alpha x,\lambda) + \zeta_2(y,\lambda) +
\zeta_2(\lambda y,\alpha).
$$

\begin{proposition}\label{prop:zeta12}
The functions $\zeta_1(x,\lambda)$ and $\zeta_2(y,\lambda)$ are
equivalently defined by
\begin{eqnarray}\label{eq:zeta12zeta12}
\zeta_1(x,\alpha) &=& \zeta_1(x) - \zeta_1(\alpha x),\nonumber\\
\zeta_2(x,\alpha) &=& \zeta_2(x) - \zeta_2(\alpha x),
\end{eqnarray}
with a slight abuse of notation.
\end{proposition}

\begin{proof}
The first property of the functions $\zeta_1(x,\lambda),\
\zeta_2(y,\lambda)$ implies
$$\zeta_1(x,1) =- \zeta_2(y,1) =\zeta\in\R$$ which is set to 0
without loss of generality. The second property leads to
\begin{eqnarray}\label{eq:eta1eta2}
\zeta_1(x,\alpha\lambda) &=& \zeta_1(x,\alpha)
+ \zeta_1(\alpha x,\lambda) + F(\alpha,\lambda),\nonumber\\
\zeta_2(y,\alpha\lambda) &=& \zeta_2(y,\lambda) + \zeta_2(\lambda
y,\alpha) - F(\alpha,\lambda),
\end{eqnarray}
for some function $F:\R^2\to\R$. Moreover, by deriving the latter
expression involving $\zeta_1(x,\lambda)$ successively by $x_i$ and
$\lambda$, we get
\begin{equation}\label{eq:eta1dd}
\partI{}{x_i}\partI{}{\lambda} \zeta_1(\alpha x,\lambda) =
\partI{}{x_i}\partI{}{\lambda} \zeta_1(x,\alpha\lambda).
\end{equation}
On the other hand, taking the derivatives with respect to $\lambda$
and $\alpha$ leads to
$$
\zeta_1'(x,\alpha\lambda) + \alpha\lambda\zeta_1''(x,\alpha\lambda)
= x\cdot\nabla_{\alpha x}\zeta_1'(\alpha x,\lambda) +
\partI{}{\alpha}\partI{}{\lambda} F(\alpha,\lambda),
$$
with $'$ denoting the derivative with respect to the last coordinate
in $\zeta_1(x,\lambda)$. Using (\ref{eq:eta1dd}), we obtain the
following:
$$
\partI{}{\alpha}\partI{}{\lambda} F(\alpha,\lambda) =
\zeta_1'(x,\alpha\lambda) + \alpha\lambda\
\zeta_1''(x,\alpha\lambda) -
x\cdot\nabla_x\zeta_1'(x,\alpha\lambda),
$$
which shows that the second derivative
$\partI{}{\alpha}\partI{}{\lambda}F(\alpha, \lambda)$ is a function
of the product $\alpha\lambda$ only. Hence with a slight abuse of
notation,
$$
F(\alpha,\lambda) = F(\alpha\lambda) - F_1(\alpha) - F_2(\lambda),
$$
where we set $F_1(1) = F_2(1) =0$ without loss of generality.
Setting respectively $\alpha=1$ and $\lambda=1$ yields from
(\ref{eq:eta1eta2}), $F_1(\lambda) = F(\lambda)$ and $F_2(\lambda) =
F(\lambda)$, hence $F(1)=0$. Equation (\ref{eq:eta1eta2}) then
becomes
\begin{eqnarray*}
\zeta_1(x,\alpha\lambda) - F(\alpha\lambda) &=& \zeta_1(x,\alpha) -
F(\alpha) + \zeta_1(\alpha x,\lambda) - F(\lambda),\\
\zeta_2(y,\alpha\lambda) + F(\alpha\lambda) &=& \zeta_2(y,\lambda) +
F(\lambda) + \zeta_2(\lambda y,\alpha) + F(\alpha).
\end{eqnarray*}
$F(\lambda)$ turns out to simply shift the functions
$\zeta_1(x,\lambda)$ and $\zeta_2(y,\lambda)$. It can thus be
included as part of those functions, so without loss of generality
$F(\lambda)=0$ and the latter equations become
\begin{eqnarray}
\zeta_1(x,\alpha\lambda) &=& \zeta_1(x,\alpha)
+ \zeta_1(\alpha x,\lambda),\nonumber\\
\zeta_2(y,\alpha\lambda) &=& \zeta_2(y,\lambda) + \zeta_2(\lambda
y,\alpha).
\end{eqnarray}
The latter functions can thus be further reduced to functions of one
variable. With a slight abuse of notation, setting $\zeta_i(x) =
\zeta_i\(x,\frac1{||x||}\)$ and choosing respectively $\lambda =
\frac1{\alpha||x||}$ and $\alpha = \frac1{\lambda||y||}$ yields the
resulting equations (\ref{eq:zeta12zeta12}).
\end{proof}

This last proposition combined with (\ref{eq:kscaling1}) gives
$$
k(\lambda x,y)+\zeta_1(\lambda x)-\zeta_2(y) = k(x,\lambda
y)+\zeta_1(x)-\zeta_2(\lambda y).
$$
Now setting $\Phi\(x,y\) = \exp
\big(k(x,y)+\zeta_1(x)-\zeta_2(y)\big)$ completes the proof of the
first condition of the theorem, i.e.
$$
\Phi(\lambda x, y) = \Phi(x,\lambda y).
$$

\subsubsection{The function $\rho(x)$}

Going back to (\ref{eq:ldt}), we replace the explicit form of
$l(x,y,t)$ and use the scaling property of $\zeta_1(x,\lambda),\
\zeta_2(y,\lambda)$ and  $k(x,y)$ to obtain
\begin{eqnarray}\label{eq:etadt}
&\displaystyle \partI{}{t} &\[ n\ln \frac1t + \zeta_1(x,\sqrt
t)-\zeta_1(x,\sqrt{t-h}) - \zeta_2\(b,\frac1{\sqrt t}\) +
\zeta_2\(a,\frac1{\sqrt{t-h}}\) \right.\nonumber\\ && +\
\zeta_1\(\frac{b}{\sqrt h},\sqrt\frac{t-h}{t}\) -
\zeta_2\(\frac{a}{\sqrt h},\sqrt\frac{t}{t-h}\)\nonumber\\
&& +\ \varphi_1\(x,\frac{1}{t}\) - \varphi_1\(x,\frac{1}{t-h}\)
 +
\varphi_1\(\frac{b}{t},\frac{h}{t(t-h)}\)\nonumber\\
&& \left.+\ \varphi_2\(\frac{b}{t},\frac{1}{t}\)-
\varphi_2\(\frac{a}{t-h},\frac{1}{t-h}\)
+\varphi_2\(\frac{a}{t-h},\frac{h}{t(t-h)}\)\]=0.
\end{eqnarray}
Note that the terms composed of the function $k(x,y)$ do not depend
on $t$ and thus cancel out with the time derivative.

Recall that the variables $a,b,x$ are independent from each other.
So taking the derivative of (\ref{eq:etadt}) with respect to $x_i$
gives
$$
t_1^2\ \partI{}{x_i}\varphi_1'(x,t_1) - \frac{\sqrt{t_1}}{2}
\partI{}{x_i}\zeta_1'\(x,\frac1{\sqrt{t_1}}\) = t_2^2\
\partI{}{x_i}\varphi_1'(x,t_2)
- \frac{\sqrt{t_2}}{2}\partI{}{x_i}\zeta_1'\(x,\frac1{\sqrt{t_2}}\)
$$
with $t_1=\frac1t$, $t_2=\frac1{t-h}$ and $'$ still denotes the
derivative with respect to the last coordinate. Since the latter is
valid for all $t_1, t_2>0$, each side of the equation must be
independent of $t$. Hence, for some differentiable function
$\varphi_{11}:\R^n\to\R$,
$$
\partI{}{t}\partI{}{x_i}\varphi_1(x,t) = -\frac1{t^2}
\partI{}{x_i}\varphi_{11}(x)
- \partI{}{t}\partI{}{x_i}\zeta_1\(x,\frac1{\sqrt t}\),
$$
which integrates to
$$
\varphi_1(x,t) = \frac1t \varphi_{11}(x) - \zeta_1\(x,\frac1{\sqrt
t}\) + h_1(x) + \tau_1(t).
$$
There are so far no further conditions to add on the function
$\varphi_{1}(x,t)$.

With the new result for the shape of $\varphi_1(x,t)$ and using the
scaling property of $\zeta_1(x,\lambda)$ and $\zeta_2(y,\lambda)$,
we derive (\ref{eq:etadt}) with respect to $b_i$ to obtain
\begin{eqnarray*}
&\displaystyle
\partI{}{b_i}\partI{}{t} & \[\zeta_1\(b,\frac1t\)
-\zeta_2\(b,\frac1{\sqrt t} \) +
\frac{t(t-h)}{h}\varphi_{11}\(\frac{b}{t}\)
\right.\\
&&\left. +\  h_1\(\frac{b}{t}\) +
\varphi_2\(\frac{b}{t},\frac{1}{t}\)\] =0.
\end{eqnarray*}
For convenience, let $\varphi_{22}:\R^{n+1}\to\R$ be such that
\begin{equation}
\varphi_2(x,t) = \zeta_2\(\frac{x}{t},\sqrt
t\)-\zeta_1\(\frac{x}{t},t\) + \frac1t\varphi_{11}(x)- h_1(x)+
\varphi_{22}(x,t).
\end{equation}
The former partial derivative equation then becomes
\begin{equation}\label{eq:phi22b}
\partI{}{b_i}\partI{}{t} \[
\frac{t^2}{h}\varphi_{11}\(\frac{b}{t}\) +
\varphi_{22}\(\frac{b}{t},\frac{1}{t}\)\] =0,
\end{equation}
which develops to
$$
t_1^2\left\{1 + t_1\partI{}{t_1} + z\cdot\nabla\right\}
\partial_i\varphi_{22}(z,t_1) = \frac1h\bigg\{1 -
z\cdot\nabla\bigg\} \partial_i\varphi_{11}(z),
$$
with the notation $z=\frac{b}{t},\ t_1=\frac{1}{t}$. Since the LHS
is independent of $h$, both sides of the equation must cancel out.
This is the case if $\partial_i\varphi_{22}(z,t)$ is homogeneous of
degree $-1$ and $\partial_i\varphi_{11}(z)$ is homogeneous of degree
1, which means by integration over $z_i$,
\begin{equation}
\varphi_{11}(\lambda z) = \lambda^2\varphi_{11}(z) +
\bar\varphi\ln\lambda, \quad \bar\varphi\in\R.
\end{equation}
Let $\rho(z) = \varphi_{11}(z)$, then, conditioned on showing
$\bar\varphi=0$, we recover the second condition of the theorem,
that is
$$
\rho(\lambda z) = \lambda^2\rho(z).
$$

\subsubsection{The function $\theta(y)$}

In order to further investigate the properties of
$\varphi_{22}(z,t)$, we derive (\ref{eq:etadt}) by $a_i$:
\begin{eqnarray*}
&\displaystyle
\partI{}{a_i}\partI{}{t} & \[
\varphi_{22}\(\frac{a}{t-h},\frac{h}{t(t-h)}\) -
\varphi_{22}\(\frac{a}{t-h},\frac{1}{t-h}\)\right.\\
&&\left. +\ \zeta_1\(a,\frac{t}{h}\) - \zeta_2\(a,\frac{t}{h}\)\]
=0,
\end{eqnarray*}
using the scaling properties of $\zeta_1(x,\lambda),\
\zeta_2(y,\lambda)$ and $\varphi_{11}(x)$. Since
$\partial_i\varphi_{22}(z,t)$ is homogeneous of degree $-1$, this
leads to
$$
\partI{}{t}\[\partial_i\varphi_{22}\(a,\frac{h}{t}\)
+\partial_i\zeta_1\(a,\frac{t}{h}\) -
\partial_i\zeta_2\(a,\frac{t}{h}\)\]=0.
$$
Hence, by integration over time, we obtain
$$
\partial_i\varphi_{22}\(a,\frac{h}{t}\)
+\partial_i\zeta_1\(a,\frac{t}{h}\) -
\partial_i\zeta_2\(a,\frac{t}{h}\)=\partial_i h_2(a),
$$
for some differentiable function $h_2:\R^n\to\R$. Integration with
respect to $a_i$ yields
$$
\varphi_{22}\(a,t\) = h_2(a) + \tau_2(t)-\zeta_1\(a,\frac1t\) +
\zeta_2\(a,\frac1t\),
$$
for some real function $\tau_2(t)$. Furthermore the homogeneity
condition on $\partial_i\varphi_{22}(z,t)$ translates to $\partial_i
h_2(z)$ as follows:
$$
\lambda\partial_i h_2(\lambda a) = \partial_i h_2(a) + \partial_i
\zeta_2(a,\lambda) - \partial_i \zeta_1(a,\lambda),
$$
by the scaling properties of $\partial_i\zeta_1(x,\lambda)$ and
$\partial_i \zeta_2(y,\lambda)$. Integration over $a_i$ then yields
$$
h_2(\lambda a) = h_2(a) + H(\lambda) + \zeta_2(a,\lambda) -
\zeta_1(a,\lambda),
$$
for $H:\R\to\R$ with the properties $H(1) =0$ and
$H(\lambda_1\lambda_2) = H(\lambda_1)+H(\lambda_2)$. The function
$H(\lambda)$ is thus of the form $H(\lambda) = \beta\ln\lambda$ for
some $\beta\in\R$. Using Proposition \ref{prop:zeta12}, the latter
becomes
$$
h_2(\lambda a) + \zeta_2(\lambda a) - \zeta_1(\lambda a) =
\beta\ln\lambda + h_2(a) + \zeta_2(a) - \zeta_1(a).
$$
Let $\theta(z) = \exp\big( h_2(z)+ \zeta_2(a) - \zeta_1(a) \big)$.
We then recover condition 3 of the theorem, i.e.
$$
\theta(\lambda z) = \lambda^{\beta}\ \theta(z).
$$

\subsubsection{The function $\tau(t)$}
The results for $\varphi_1(x,t)$ and $\varphi_2(y,t)$ summarize so
far to:
\begin{eqnarray}
\varphi_1(x,t) &=& \frac{\rho(x)}{t} - \zeta_1(x) +
\zeta_1\(\frac{x}{\sqrt t}\) + h_1(x) + \tau_1(t)\nonumber\\
\varphi_2(y,t) &=& \frac{\rho(y)}{t} +\zeta_2(y) -
\zeta_2\(\frac{y}{\sqrt t}\) - h_1(y) + \tau_2(t) + h_2(y).\nonumber
\end{eqnarray}
This implies that $l(x,y,t) = \ln( p_t(x,y) )$ has the form
\begin{equation}
l(x,y,t) = \ln\Phi\(\frac{x}{\sqrt t},\frac{y}{\sqrt t}\) +
\ln\frac{\eta(y)}{\eta(x)} + \frac{\rho(x) + \rho(y)}{t} + \tau_1(t)
+ \tau_2(t),
\end{equation}
where $\eta(x)= \exp\big(\zeta_1(x) - h_1(x)\big)$ defines an
$h$-transform and can thus be neglected. It remains to explicitly
formulate the form of the function $\tau:\R_+\to\R$ that we define
by
$$
\tau(t) = \tau_1(t) + \tau_2(t) +\frac12(n+\beta)\ln t
$$
to satisfy the resulting partial differential equation
(\ref{eq:etadt}):
\begin{eqnarray}
&\displaystyle \partI{}{t} &\[
\bar\varphi\ \frac{t(t-h)}{h}\(\frac{t}{t-h}\ln \frac1t
+ \frac{t-h}{t}\ln\frac{1}{t-h}\) \right.\nonumber\\
&& \left.+\ \tau\(\frac1t\) +
\tau\(\frac{h}{t(t-h)}\)-\tau\(\frac{1}{t-h}\)\]=0.
\end{eqnarray}
The former equation becomes after derivation,
$$
\frac{\bar\varphi}{t_2-t_1}\Big(t_1(2\ln t_2+1) + t_2(2\ln t_1
+1)\Big) = t_1^2\tau'(t_1) +(t_2^2 - t_1^2) \tau'(t_2-t_1) -
t_2^2\tau'(t_2),
$$
for $t_1=\frac1t$, $t_2=\frac1{t-h}$, which hints that
$\bar\varphi=0$. Assuming $\bar\varphi$ to be non-zero and setting
$\displaystyle g(t) = \frac{t^2}{2\bar\varphi}\ \tau'(t)+\frac12$,
we get indeed
$$
t_1\ln t_2 + t_2\ln t_1 = t_1 \Big(g(t_2-t_1)- g(t_1) + g(t_2)\Big)
+ t_2 \Big(g(t_2-t_1)+ g(t_1) - g(t_2)\Big),
$$
which has no solution. So $\bar\varphi=0$. Now dividing the former
equation with $\bar\varphi=0$ by $t_1$ and taking the limit
$t_1\to0$ leads to the differential equation
$$
\tau''(t) = 0,
$$
which solves for $\tau(t) = \nu t,\ \nu\in\R$. However, $e^{\nu t}$
can always be included in an $h$-transform, so $\nu$ is set to 0 and
$\tau(t)$ becomes trivial.

\subsubsection{The general case $\alpha>0$}
Combining the different factors, we obtain for the semigroup
densities,
$$
p_t(x,y) = t^{-\frac{1}{2}\(n+\beta\)}\ \Phi\(\frac{x}{\sqrt t},
\frac{y}{\sqrt t}\)\ \theta(y)\ \exp\bigg\{\frac{\rho(x)}{t} +
\frac{\rho(y)}{t}\bigg\},
$$
which simplifies to, using the scaling property of $\theta(y)$ and
$\rho(x)$,
$$
p_t(x,y) = t^{-\frac{n}{2}}\ \Phi\(\frac{x}{\sqrt t}, \frac{y}{\sqrt
t}\)\ \theta\(\frac{y}{\sqrt t}\)\ \exp\bigg\{\rho\(\frac{x}{\sqrt
t}\) +\rho\(\frac{y}{\sqrt t}\)\bigg\}.
$$
For $\alpha>0$, we apply the change of variables $x\mapsto
\|x\|^\alpha \frac{x}{\|x\|}$. The semigroup densities become
$$
p_t^\alpha(x,y) = p_t\(x^\frac{1}{\alpha}, y^\frac{1}{\alpha}\) J(y)
$$
where $J(y)$ is the Jacobian of the inverse transformation, that is
$J(y) = \alpha^{-n}\ y_1^{\frac{1}{\alpha}-1}\cdots
y_n^{\frac{1}{\alpha}-1}$, and we use the slight abuse of notation
$x^\frac{1}{\alpha} = \(x_1^\frac{1}{\alpha},
\ldots,x_n^\frac{1}{\alpha}\)$. The semigroup densities can be
recast into
$$
p_t^\alpha(x,y) = t^{-\frac{n\alpha}{2}}\
\tilde\Phi\(\frac{x}{t^\frac{\alpha}{2}},
\frac{y}{t^\frac{\alpha}{2}}\)\
\tilde\theta\(\frac{y}{t^\frac{\alpha}{2}}\)\
\exp\bigg\{\tilde\rho\(\frac{x}{t^\frac{\alpha}{2}}\)
+\tilde\rho\(\frac{y}{t^\frac{\alpha}{2}}\)\bigg\}
$$
where $\tilde\Phi(x,y) = \Phi\(x^\frac{1}{\alpha},
y^\frac{1}{\alpha}\)$ satisfies condition 1, $\tilde\theta(y) =
\theta\(y^\frac{1}{\alpha}\) J(y)$ satisfies condition 2 for
$\tilde\beta = \frac{\beta+n}{\alpha}-n$ and $\tilde\rho(x) =
\rho\(y^\frac{1}{\alpha}\)$ satisfies condition 3 of equation
(\ref{eq:density}). 

\section{Application to diffusions on $\R_+$}
\label{sec:applications}

The case of the diffusions on $\R_+$ was entirely characterized by
Watanabe in \cite{Wa1975}. It was shown that only Bessel processes
in the wide sense (which we recall the definition below) enjoy the
time-inversion property of degree 1.
\begin{definition}
For some $\nu>-1$ and $c\ge0$, the diffusion process generated by
\begin{equation}\label{eq:BesselLwide}
\L = \frac{1}{2}\partII{}{x} + \(\frac{2\nu+1}{2x} +
\frac{h_c'(x)}{h_c(x)}\)\partI{}{x}
\end{equation}
is called {\it Bessel process in the wide sense}. The function
$h_c(x)$ is given by
\begin{equation}
h_c(x) = 2^\nu \Gamma(\nu+1)\ (\sqrt{2c}\ x)^{-\nu}\
I_\nu(\sqrt{2c}\ x),
\end{equation}
where $I_\nu$ is the modified Bessel function.
\end{definition}

\begin{remark}
The Bessel process in the wide sense is in $h$-transform
relationship, for $h\equiv h_c$, with the Bessel process.
\end{remark}

We show that the result in \cite{Wa1975} is a consequence of Theorem
\ref{thm:main}, which has the following one-dimensional formulation:
\begin{theorem}\label{thm:1d}
The Markov process $(t^\alpha X_\frac{1}{t}, t>0)$ on $\R_+$ is
homogeneous if and only if the semigroup densities of $(X_t,t\ge0)$
are of the form:
\begin{equation}\label{eq:density1d}
p_{t}(x,y) = t^{-\frac{\alpha}{2}(1+\beta)}\
\phi\(\frac{xy}{t^\alpha}\)\ y^\beta\
\exp\left\{-\frac{k^2}{2}\(\frac{x^\frac{2}{\alpha}}{t} +
\frac{y^\frac{2}{\alpha}}{t}\)\right\}\
\end{equation}
for $k>0$, or if it is in $h$-transform relationship with it.
Moreover, the semigroup densities are related as follows:
\begin{equation}\label{eq:qhp1d}
q_{t}^{(x)}(a,b) = \frac{\phi\(xb\)}{\phi\(xa\)}\ \exp\(-t\
\frac{k^2}{2}\ x^\frac{2}{\alpha}\)\ p_t(a,b).
\end{equation}
\end{theorem}

\begin{proof}
Theorem \ref{thm:main} formulated on $\R_+$ gives $ \Phi(x,y) =
\phi(xy)$, $\rho(x) = -\frac{k^2}{2}\ x^\frac{2}{\alpha}, \ k>0$,
and $\theta(y) = y^\beta$ for conditions 1, 2 and 3 to be satisfied.
Equation (\ref{eq:qhp1d}) is then a consequence of (\ref{eq:qp}).
\end{proof}

We identify further the class of diffusion processes and provide a
different proof for the result in \cite{Wa1975}.

\begin{proposition}\label{prop:Bessel}
If $(X_t,t\ge0)$ is a diffusion process and $(t X_\frac{1}{t}, t>0)$
is homogeneous and conservative, then both are necessarily (possibly
time-scaled) Bessel processes in the wide sense.
\end{proposition}

\begin{proof}
If $(X_t,t\ge0)$ is a diffusion process, then its infinitesimal
generator has the following general structure:
$$
\L = \frac{d}{\mathfrak{m}(dy)}\frac{d}{d\mathfrak{s}(y)}
$$
where $\mathfrak{m}(dy)$ is called the speed measure and
$\mathfrak{s}(y)$ the scale measure density (see \cite{IM1974}).
From the assumptions of the proposition, $(X_t,t\ge0)$ enjoys the
time-inversion property of degree 1. As a consequence of
(\ref{eq:density1d}) in Theorem \ref{thm:1d}, the speed measure is
thus absolutely continuous with respect to the Lebesgue measure. Let
$\mu(x)=\frac{1}{\mathfrak{m}(x)\mathfrak{s}(x)}$ and
$s(x)=\frac{-\mathfrak{s}'(x)} {\mathfrak{m}(x)\mathfrak{s}(x)^2}$.
The infinitesimal generator can then be expressed as
$$
\L = s(y)\partII{}{y} + \mu(y)\partI{}{y}
$$
where it remains to identify the functions $s(y)$ and $\mu(y)$. For
a fixed $x>0$, let $\L^{(x)}$ be the infinitesimal generator of $(t
X_\frac{1}{t}, t>0)$. From equation (\ref{eq:qhp1d}), $\L^{(x)}$ has
the following relationship with $\L$:
$$
\L^{(x)}:f(b) \mapsto \frac{1}{\phi(xb)}\ \L \big(\phi(xb) f(b)\big)
- \frac{k^2}{2}\ x^2 f(b),
$$
which develops to
$$
\L^{(x)}f(b) = s(b) f''(b) + \left\{\mu(b) + 2x\ s(b)\
\frac{\phi'(xb)}{\phi(xb)}\right\}f'(b) + U(x,b) f(b).
$$
For the process to be conservative, we require $U(x,b)=0$, which
implies no killing in the interior of the domain, that is
$$
s(b)\ x^2\ \frac{\phi''(xb)}{\phi(xb)} + \mu(b)\ x\
\frac{\phi'(xb)}{\phi(xb)} - \frac{k^2}{2}\ x^2 =0.
$$
We change variables to $z=xb$ to obtain
$$
s(z/x)\ \phi''(z) + \frac{\mu(z/x)}{x}\ \phi'(z) - \frac{k^2}{2}\
\phi(z) =0.
$$
Since the latter must be valid for all $x>0$, we are led to set:
$s(b)=\frac{\sigma^2}{2}$ for $\sigma>0$ and $\mu(b) =
\frac{\sigma^2}{2}\ \frac{2\nu+1}{b}$ for $\nu>-1$. This yields the
following equation
$$
\frac{1}{2}\phi''(z) + \frac{2\nu+1}{2z}\ \phi'(z) -
\frac{k^2}{2\sigma^2}\ \phi(z) = 0.
$$
The general solution (non-singular at 0 and up to a constant factor)
is expressed through the modified Bessel function of the first kind
as follows:
$$
\phi(z) = \(\frac{k z}{\sigma}\)^{-\nu}\ I_\nu\(\frac{k
z}{\sigma}\).
$$
Gathering the different factors in (\ref{eq:density1d}) leads to
$$
p_{t}(x,y) = N\ t^{-\frac{1+\beta}{2}}\ \(\frac{xy}{t}\)^{-\nu}\
I_\nu\(\frac{k}{\sigma}\frac{xy}{t}\)\ y^{\beta}\
\exp\left\{-\frac{k^2}{2}\(\frac{x^2}{t} + \frac{y^2}{t}\)\right\},
$$
where $N$ is a normalization factor. The additional condition that
$\displaystyle \lim_{t\to0} p_t(x,y) = \delta(x-y)$ implies $\beta =
2\nu +1$ and $k=\frac{1}{\sigma}$, which leads to the semigroup
densities of a time-scaled Bessel process of dimension $\nu$:
$$
p_{t}(x,y) = \frac{y}{\sigma^2t}\ \(\frac{y}{x}\)^{\nu} \
I_\nu\(\frac{xy}{\sigma^2t}\)\ \exp\left\{-\frac{x^2+y^2}{2\sigma^2
t}\right\}.
$$
The infinitesimal generator for $(t X_\frac{1}{t}, t>0)$ is given by
$$
\L^{(x)} = \frac{\sigma^2}{2} \partII{}{b} +
\sigma^2\left\{\frac{2\nu+1}{2b} +
x\frac{\phi'(xb)}{\phi(xb)}\right\}\partI{}{b},
$$
where one recognizes expression (\ref{eq:BesselLwide}) for $h_c(b) =
\phi(xb)$ and a time-scale $t\mapsto\sigma^2 t$.
\end{proof}

Proposition \ref{prop:Bessel} smoothly extends to any power
$\alpha>0$ of the Bessel process through the mapping $x\mapsto
||x||^\alpha\frac{x}{||x||}$. In particular, it is worth remarking
that the case $\alpha=2$ gives rise to squares of Bessel processes,
which leads to the following result:
\begin{proposition}\label{prop:Bessel2}
If $(X_t,t\ge0)$ is a diffusion process and $(t^2 X_\frac{1}{t},
t>0)$ is homogeneous and conservative, then both are necessarily
(possibly time-scaled) squares of Bessel processes in the wide
sense.
\end{proposition}

\section{Examples}
\label{sec:examples}

\subsection{Generalized Dunkl processes and Jacobi-Dunkl processes}

\subsubsection{Multidimensional Dunkl processes}
We briefly review the construction of the Dunkl process in $\R^n$
(see \cite{Ro1998,RV1998}).
\begin{definition}
The Dunkl process in $\R^n$ is the Markov c{\`a}dl{\`a}g process
with infinitesimal generator
\begin{equation}
\frac{1}{2}\L^{(k)} = \frac{1}{2}\sum_{i=1}^n T_i^2
\end{equation}
where $T_i,\ 1\le i\le n,$ is a one-dimensional
differential-difference operator defined for $u\in C^1(\R^n)$ by
\begin{equation}
T_i u(x) = \partI{u(x)}{x_i} + \sum_{\alpha\in R_+}
k(\alpha)\alpha_i \frac{u(x)-u(\sigma_\alpha
x)}{\langle\alpha,x\rangle}.
\end{equation}
$\langle\cdot,\cdot\rangle$ is the usual scalar product. $R$ is a
root system in $\R^n$ and $R_+$ a positive subsystem. $k(\alpha)$ is
a non-negative multiplicity function defined on $R$ and invariant by
the finite reflection group $W$ associated with $R$. $\sigma_\alpha$
is the reflection operator with respect to the hyperplane $H_\alpha$
orthogonal to $\alpha$ such that $\sigma_\alpha x = x -
\langle\alpha,x\rangle\alpha$ and for convenience
$\langle\alpha,\alpha\rangle=2$ (see \cite{Du1989,Du1992}).
\end{definition}

A result obtained by M. R\"osler \cite{Ro1998} yields the semigroup
densities as follows:
\begin{equation}\label{eq:DunklKernel}
p^{(k)}_t(x,y) = \frac{1}{c_k t^{\gamma+n/2}}\
\exp\(-\frac{|x|^2+|y|^2}{2t}\)\ D_k\(\frac{x}{\sqrt t},
\frac{y}{\sqrt t}\)\ \omega_k(y)
\end{equation}
where $D_k(x,y)>0$ is the Dunkl kernel, $\displaystyle\omega_k(y) =
\prod_{\alpha\in R_+}|\langle\alpha,y\rangle|^{2k(\alpha)}$ the
weight function which is homogeneous of degree $\displaystyle
2\gamma = 2\sum_{\alpha\in R_+} k(\alpha)$ and $\displaystyle c_k =
\int_{\R^n} e^{-\frac{|x|^2}{2}}\omega_k(x)dx$.

Following a thorough study of the properties of the one-dimensional
Dunkl process in \cite{GY2005a}, it was remarked in \cite{GY2005b}
that the Dunkl process in $\R^n$ enjoys the time-inversion property
of degree 1. Considering that the Dunkl kernel satisfies
\begin{equation}\label{eq:DunklKerProp}
D_k(x,y) = D_k(y,x)\quad\mathrm{and}\quad D_k(\mu x,y) = D_k(x,\mu
y),
\end{equation}
the proof is straightforward with
\begin{equation}
\Phi(x,y) \equiv D_k(x,y),\quad \theta(y) \equiv \omega_k(y),
\quad\rho(x)\equiv -\frac{|x|^2}{2}.
\end{equation}
By equation (\ref{eq:qhp}), the semigroup densities of the
time-inverted process is even in $h$-transform relationship with the
semigroup densities of the original Dunkl process:
\begin{equation}
q_t^{(x)}(a,b) = \frac{D_k(x,b)}{D_k(x,a)}\ \exp\(-\frac{|x|^2}{2}
t\)\ p_t^{(k)}(a,b).
\end{equation}

\subsubsection{Generalized Dunkl processes}

In an attempt to generalize the Dunkl process, we extend the
definition of the infinitesimal generator to
\begin{equation}\label{eq:GenDunklGen}
\L^{(k,\lambda)} f(x) = \frac{1}{2}\Delta f(x) + \sum_{\alpha\in
R_+} k(\alpha)\frac{\langle\nabla f(x),\alpha\rangle}{\langle x
,\alpha\rangle} + \sum_{\alpha\in R_+} \lambda(\alpha)\frac{
f(\sigma_\alpha x) - f(x)}{\langle x,\alpha\rangle^2}
\end{equation}
where $\Delta$ is the usual Laplacian, $f\in C^2(\R^n)$ and
$\lambda(\alpha)$ is a non-negative multiplicity function defined on
$R$ and invariant by the finite reflection group $W$, similarly as
$k(\alpha)$. We retrieve the Dunkl process for $\lambda(\alpha) =
k(\alpha)$. Note that these processes are no longer martingales for
$\lambda(\alpha) \neq k(\alpha)$.

The one-dimensional case was introduced in \cite{GY2005a}, where the
semigroup densities were explicitly derived,
\begin{equation}
p_t^{(k,\lambda)}(x,y) = \frac{1}{t^{k-\frac{1}{2}}}\
y^{k-\frac{1}{2}} \exp\(-\frac{x^2+y^2}{2 t}\)\
D_{k,\lambda}\(\frac{xy}{t}\)
\end{equation}
with the generalized Dunkl kernel, $(\nu=k-\frac{1}{2},\ \mu =
\sqrt{\nu^2+4\lambda})$,
\begin{equation}
D_{k,\lambda}\(z\) =
1_{\{y\in\R_-\}}\frac{1}{2z^\nu}(I_\nu-I_\mu)\(-z\) +
1_{\{y\in\R_+\}}\frac{1}{2z^\nu}(I_\nu+I_\mu)\(z\),
\end{equation}
for $z=\frac{xy}{t}$. From (\ref{eq:DunklKerProp}), the generalized
Dunkl kernel satisfies
\begin{equation}
D_{k,\lambda}(x,y) = D_{k,\lambda}(y,x)\quad\mathrm{and}\quad
D_{k,\lambda}(\mu x,y) = D_{k,\lambda}(x,\mu y),
\end{equation}
which readily implies that the generalized Dunkl process also enjoys
the time-inversion property of degree 1. The semigroup densities
were derived as an application of the skew-product representation of
the generalized Dunkl process $(X_t,t\ge0)$ in terms of its absolute
value (a Bessel process) and an independent Poisson process
$N_t^{(\lambda)}$:
\begin{equation}
X_t \stackrel{(d)}{=}  |X_t| (-1)^{N^{(\lambda)}_{A_t}}
\end{equation}
where $A_t = \int_0^t \frac{ds}{X^2_s}$.

In the $n$-dimensional case, the application of the skew-product
representation derived by Chybiryakov \cite{Ch2005} shows that the
generalized Dunkl process enjoys time-inversion for some specific
root systems. We first recall one of the main results of
\cite{Ch2005}.

\begin{proposition}\label{prop:DunklSkewProd}
Let $(X_t,t\ge0)$ be the generalized Dunkl process generated by
(\ref{eq:GenDunklGen}), with $(X_t^W,t\ge0)$ its radial part, i.e.
the process confined to a Weyl chamber. Let
$R_+\equiv\{\alpha_1,\ldots,\alpha_l\}$ for some $l\in\N$ be the
corresponding positive root system and let $(N_t^i,t\ge0),\
i=1,\ldots,l$ be independent Poisson processes of respective
intensities $\lambda(\alpha_i)$. Then $X_t$ may be represented as
$Y_t^l$, which is defined by induction as follows:
$$
Y_t^0 = X_t^W\ \ {\rm and}\ \ Y_t^i =
\sigma_{\alpha_i}^{N^i_{A_t^i}}\ Y_t^{i-1},\ \ i=1,\ldots,l,
$$
where $\displaystyle A_t^i = \int_0^t \frac{ds}{\langle
Y_s^{i-1},\alpha_i\rangle^2}$.
\end{proposition}

\noindent The proof follows the argument in \cite{Ch2005}, while
replacing $k(\alpha_i)$ by $\lambda(\alpha_i)$ appropriately.

From now on, let $R_+\equiv\{\alpha_1,\ldots,\alpha_l\}$ for some
$l\le n$ be an orthogonal positive root system, that is
$\langle\alpha_i,\alpha_j\rangle=2\delta_{ij}$. For this particular
root system, one can show that the generalized Dunkl process enjoys
time-inversion of degree 1. We first prove the following absolute
continuity relation:

\begin{lemma}\label{lem:DunklExpec}
Let $(X_t^W,t\ge0)$ be the radial Dunkl process with infinitesimal
generator
\begin{equation}
\L^W_k f(x) = \frac{1}{2}\Delta f(x) + \sum_{i=1}^l
k(\alpha_i)\frac{\langle\alpha_i,\nabla f(x)\rangle}
{\langle\alpha_i,x\rangle}.
\end{equation}
Fixed $\nu\in\{1,\ldots,l\}$. Let $k'(\alpha)$ be another
coefficient function on the root system $R_+$ such that
$k'(\alpha_\nu)>k(\alpha_\nu)$ and $k'(\alpha_i)=k(\alpha_i)$ for
$i\neq\nu$. Then, denoting $\P_x^{(k)}$ the law of the radial Dunkl
process $X_t^W$ starting from $x$, we have
\begin{eqnarray}\label{eq:DunklAbsCont}
\P^{(k')}_x\big|_{\F_t} &=& \(\frac{\langle\alpha_\nu,X_t^W\rangle}
{\langle\alpha_\nu,x\rangle}\)^{k'(\alpha_\nu)-k(\alpha_\nu)}\\
&\cdot&\exp\[- \frac{\(k'(\alpha_\nu)-\frac{1}{2}\)^2-
\(k(\alpha_\nu)-\frac{1}{2}\)^2}{2}\int_0^t \frac{ds}{\langle
\alpha_\nu,X_s^W\rangle^2}\]\cdot \P^{(k)}_x\big|_{\F_t}.\nonumber
\end{eqnarray}
\end{lemma}

\begin{proof}
Let $k_0(\alpha)$ be a coefficient such that $k_0(\alpha_\nu) =
\frac{1}{2}$ for some fixed $i\in\{1,\ldots,l\}$. $X_t^W$ has the
following martingale decomposition (see \cite{GY2004c}):
$$
X_t^W = x + B^{(k_0)}_t + \sum_{i=1}^l k_0(\alpha_i)\int_0^t
\frac{ds}{\langle \alpha_i,X_s^W\rangle}\ \alpha_i
$$
where $B_t^{(k_0)}$ is a $(\P^{(k_0)}_x,\F_t)$-Brownian motion.
Consider the local martingale
$$
L^{(k')}_t = \exp\(\(k'(\alpha_\nu)-\frac{1}{2}\) \int_0^t
\frac{\langle \alpha_\nu,dB_s^{(k_0)}\rangle}{\langle
\alpha_\nu,X_s^W\rangle} -
\frac{\(k'(\alpha_\nu)-\frac{1}{2}\)^2}{2}\int_0^t \frac{ds}{\langle
\alpha_\nu,X_s^W\rangle^2}\)
$$
for some coefficient function $k'(\alpha)$ such that
$k'(\alpha_\nu)>\frac{1}{2}$ and $k'(\alpha_i)=k_0(\alpha_i)$ for
$i\neq\nu$. The It\^o formula for
$\ln\(\langle\alpha_\nu,X_t^W\rangle\)$ combined with the
orthogonality of the roots yields
$$
L^{(k')}_t = \(\frac{\langle\alpha_\nu,X_t^W\rangle}
{\langle\alpha_\nu,x\rangle}\)^{k'(\alpha_\nu)-\frac{1}{2}}\ \exp\(-
\frac{\(k'(\alpha_\nu)-\frac{1}{2}\)^2}{2}\int_0^t \frac{ds}{\langle
\alpha_\nu,X_s^W\rangle^2}\).
$$
Define the new law $\P^{(k')}_x|_{\F_t} = L^{(k')}_t\cdot
\P^{(k_0)}_x|_{\F_t}$. By the Girsanov theorem,
$$
B_t^{(k')} = B_t^{(k_0)} - \(k'(\alpha_\nu)-\frac{1}{2}\)
\int_0^t\frac{ds}{\langle\alpha_\nu, X_s^W\rangle}
$$
is a $(\P^{(k')}_x,\F_t)$-Brownian motion and hence,
$$
X_t^W = x + B_t^{(k')} + \sum_{i=1}^l k'(\alpha_i) \int_0^t
\frac{ds}{\langle \alpha_i,X_s^W\rangle}\ \alpha_i
$$
is a radial Dunkl process under $(\P^{(k')}_x,\F_t)$.

Define $k$ as another coefficient on the root system that satisfies
the conditions enunciated in the lemma. The absolute continuity
relation is then a consequence of the successive application of the
latter result to the indices $k'$ and $k$.
\end{proof}

Now as an application of Proposition \ref{prop:DunklSkewProd}, we
prove the following:

\begin{proposition}
Let $R_+\equiv\{\alpha_1,\ldots,\alpha_l\}$ be an orthogonal
positive root system for $l\le n$. Then the generalized Dunkl
process $(X_t,t\ge0)$ generated by (\ref{eq:GenDunklGen}) enjoys the
time-inversion property of degree 1.
\end{proposition}

\begin{proof}
Using orthogonality of the roots, remark that
$$
\langle\alpha_i,\sigma_j x\rangle^2 =
\langle\alpha_i,x-\langle\alpha_j, x\rangle\alpha_j\rangle^2 =
\langle\alpha_i,x\rangle^2,
$$
which implies in particular
$$
\langle\alpha_i,Y_t^i\rangle^2 = \langle\alpha_i,X_t^W\rangle^2,
$$
so that the inductive representation of $X_t$ in Proposition
\ref{prop:DunklSkewProd} becomes
$$
X_t = \prod_{i=1}^l \sigma_{\alpha_i}^{N^i_{A_t^i}}\ X_t^W\quad {\rm
for}\quad A_t^i = \int_0^t \frac{ds}{\langle
X_s^W,\alpha_i\rangle^2}.
$$
The radial part of a Dunkl process enjoys the time-inversion
property of degree 1. We need to show that the semigroup densities
of the generalized Dunkl process are related to the semigroup
densities of its radial parts. For $f\in C^2(\R^n)$,
$$
\E_x\big[f(Y_t^i)\big] = \E_x\big[f(Y_t^{i-1})1_{\{N^i_{A_t^i}{\rm\
is\ even}\}}\big] + \E_x\big[f(\sigma_{\alpha_i}
Y_t^{i-1})1_{\{N^i_{A_t^i}{\rm\ is\ odd}\}}\big].
$$
Since $\P(N^i_u{\rm\ is\ even}) =
\frac{1}{2}\(1+\exp(-2\lambda(\alpha_i)u)\)$, we obtain
\begin{eqnarray*}
\E_x\big[f(Y_t^i)\big] &=& \E_x\[f(Y_t^{i-1})\
\frac{1}{2}\(1+\exp(-2\lambda(\alpha_i)A_t^i)\)\] \\
&+& \E_x\[f(\sigma_{\alpha_i} Y_t^{i-1})\
\frac{1}{2}\(1-\exp(-2\lambda(\alpha_i)A_t^i)\)\].
\end{eqnarray*}
The expectation $\E_x\big[f(X_t)\big]$ can thus be evaluated by
induction on $i\in\{1,\ldots,l\}$. It follows that the semigroup
densities of $X_t$ can be expressed as the product of the semigroup
densities of its radial parts times a function involving
expectations of the form
$$
\E_x^{(k)}\[\exp\(-2\lambda(\alpha_\nu) A_t^\nu\)\bigg| X_t^W = y\],
$$
for $\nu\in\{1,\ldots,l\}$. From Lemma \ref{lem:DunklExpec},
$$
\E_x^{(k)}\[\exp\(-2\lambda(\alpha_\nu) A_t^\nu\)\bigg| X_t^W = y\]
= \frac{p^{(k')}_t(x,y)}{p^{(k)}_t(x,y)}\
\(\frac{\langle\alpha_\nu,y\rangle} {\langle\alpha_\nu,
x\rangle}\)^{k(\alpha_\nu)-k'(\alpha_\nu)},
$$
where $k'(\alpha_i) = k(\alpha_i)$ for $i\neq\nu$ and
$k'(\alpha_\nu) = \frac{1}{2}+
\sqrt{\(k(\alpha_\nu)-\frac{1}{2}\)^2+4\lambda(\alpha_\nu)}$. The
form of the semigroup densities in (\ref{eq:DunklKernel}) implies
that the expectation is a ratio of Dunkl kernels,
$$
\E_x^{(k)}\[e^{-2\lambda(\alpha_\nu) A_t^\nu}\bigg| X_t^W = y\] =
\frac{c_k}{c_{k'}}\ \frac{D_{k'}\(\frac{x}{\sqrt t},\frac{y}{\sqrt
t}\)}{D_k\(\frac{x}{\sqrt t},\frac{y}{\sqrt t}\)}\
\frac{w_{k'}\(\frac{y}{\sqrt t}\)}{w_k\(\frac{y}{\sqrt t}\)}
\(\frac{\langle\alpha_\nu,\frac{y}{\sqrt
t}\rangle}{\langle\alpha_\nu, \frac{x}{\sqrt
t}\rangle}\)^{(k-k')(\alpha_\nu)},
$$
which reduces to
$$
\E_x^{(k)}\[e^{-2\lambda(\alpha_\nu) A_t^\nu}\bigg| X_t^W = y\] =
\frac{c_k}{c_{k'}}\ \frac{D_{k'}\(\frac{x}{\sqrt t},\frac{y}{\sqrt
t}\)}{D_k\(\frac{x}{\sqrt t},\frac{y}{\sqrt t}\)}\
\(\langle\alpha_\nu,\frac{y}{\sqrt t}\rangle \langle\alpha_\nu,
\frac{x}{\sqrt t}\rangle\)^{(k'-k)(\alpha_\nu)}
$$
by definition of $k'(\alpha)$. The conditional expectation thus
satisfies condition 1 of Theorem \ref{thm:main}. As a consequence,
the conditions of Theorem \ref{thm:main} are satisfied for
$$
\Phi(x,y) \equiv D_{k,\lambda}(x,y), \quad\theta(y) \equiv
\omega_k(y), \quad\rho(x)\equiv -\frac{|x|^2}{2}
$$
where $D_{k,\lambda}(x,y)$ is a generalized Dunkl kernel given
explicitly in terms of radial Dunkl kernels and satisfying
equivalent conditions (see (\ref{eq:DunklKerProp})).
\end{proof}

\subsubsection{Jacobi-Dunkl processes}

Gallardo et al. \cite{CGM2005} derived the Jacobi-Dunkl process as
the hyperbolic analog of the one-dimensional Dunkl process. It is
defined as the process generated by
\begin{equation}
\L^{(\alpha,\beta)} f(x) = \partII{f(x)}{x} +
\frac{A'(x)}{A(x)}\partI{f(x)}{x}
+\partI{}{x}\(\frac{A'(x)}{A(x)}\)\(\frac{f(x)- f(-x)}{2}\),
\end{equation}
where $A(x) = \(\sinh^2(x)\)^{\alpha+\frac{1}{2}}
\(\cosh^2(x)\)^{\beta+\frac{1}{2}}$. From the expression of the
semigroup densities developed in \cite{CGM2005}, this process does
not enjoy the time-inversion property. Its radial part however
corresponds to the Jacobi process of index $(\alpha,\beta)$ on
$\R_+$ (see \cite{G1997,G1998}). The infinitesimal generator of the
Jacobi process, expressed by
\begin{equation}
\L^{(\alpha,\beta)} f(x) = \frac{1}{2}\partII{f(x)}{x} +
\frac{A'(x)}{A(x)}\partI{f(x)}{x},
\end{equation}
is in $h$-transform relationship with the Laplacian operator for
$h(x) = \sqrt{A(x)}$. Since the Brownian motion enjoys the
time-inversion property of degree 1, so does the Jacobi process by
Theorem \ref{thm:main}.

\subsection{Matrix-valued processes}

\subsubsection{Eigenvalue processes}
Dyson in \cite{Dy1962} described the eigenvalues of a Hermitian
Brownian motion as the joint evolution of independent Brownian
motions conditioned never to collide (see also \cite{Gr1999} and
\cite{CL2001}). It was further remarked that the process version of
the Gaussian orthogonal ensemble does not admit such a
representation for its eigenvalues. This work was extended by
K\"onig and O'Connell in \cite{KO2001} to the process version of the
Laguerre ensemble, denominated the Laguerre process and defined as
follows:

\begin{definition}
Let $B_t$ be an $n\times m$ matrix with independent standard complex
Brownian entries. The Laguerre process is the matrix-valued process
defined by $\{X_t = B_t' B_t, t\ge0\}$, where $B_t'$ is the
transpose of $B_t$.
\end{definition}

\noindent From \cite{KO2001}, the eigenvalues of the Laguerre
process evolve like $m$ independent squared Bessel processes
conditioned never to collide. No such representation however exists
for the case where the entries of $B_t$ are real Brownian motions,
i.e. the Wishart process considered by Bru in \cite{Br1991}.

The main result of \cite{KO2001} is that both of the above mentioned
eigenvalues processes (in the complex Brownian case) can be obtained
as the $h$-transform of processes with $m$ independent components.
The joint eigenvalues process is thus in $h$-transform relationship
with a process that enjoys the time inversion property of degree 1
in the case of the $m$-dimensional Brownian motion and degree 2 in
the case of the $m$-dimensional squared Bessel process, as made
explicit in the following proposition:

\begin{proposition}
Let $p_t(x_i,y_i)\ (i=1,\ldots,m)$ be the semigroup densities of
squared Bessel processes (respectively Brownian motions), and let
\begin{equation}
h(x) = \prod_{i<j}^m (x_j-x_i)
\end{equation}
for $x=(x_1,\ldots,x_m)$. Then, the semigroup densities of the joint
eigenvalues process of the Laguerre process (respectively the
Hermitian Brownian motion) are given by
\begin{equation}
\tilde p_t(x,y) = \frac{h(y)}{h(x)} \det\big( p_t(x_i,y_j)
\big)_{i,j=1}^m
\end{equation}
with respect to the Lebesgue measure $\displaystyle dy =
\prod_{j=1}^m dy_j$.
\end{proposition}

\noindent It follows immediately by Theorem \ref{thm:main} that the
eigenvalues processes enjoy the time-inversion property. Moreover by
Proposition \ref{prop:h}, they yield the same process under
time-inversion as the $m$-dimensional Brownian motion or the
$m$-dimensional squared Bessel process respectively.

\subsubsection{Wishart processes}

The Wishart process WIS($\delta,tI_m,\frac{1}{t}x$), introduced by
Bru in \cite{Br1989a}, is a continuous Markov process taking values
in the space of real symmetric positive definite $m\times m$
matrices $S^+_m$. It is solution to the following stochastic
differential equation:
\begin{equation}
dX_t = \sqrt{X_t} dB_t + dB_t' \sqrt{X_t} + \delta I_m dt,\quad X_0
=x,
\end{equation}
where $B_t$ is an $m\times m$ matrix with Brownian entries and $I_m$
the identity matrix. Further results have been obtained in
\cite{Br1989b} and \cite{Br1991}. In \cite{DDMY2004}, among other
major findings about the Wishart process, the transition probability
densities expressed with respect to the Lebesgue measure
$\displaystyle dy = \prod_{i\le j}(dy_{ij})$ were derived in terms
of generalized Bessel functions (we refer to the appendix for the
definition):
\begin{equation}\label{eq:WishartDensities}
p_t(x,y) = \frac{1}{(2t)^{\frac{m(m+1)}{2}}}\
\exp\(-\frac{1}{2t}Tr(x+y)\)\
\(\frac{\det(y)}{\det(x)}\)^{\frac{\delta-m-1}{4}}\ {\bf\tilde
I}_{\frac{\delta-m-1}{2}} \(\frac{xy}{4t^2}\),
\end{equation}
for $x,y\in S^+_m$ and $\delta>m-1$. From the shape of its
densities, the Wishart process was stated in \cite{GY2005b} as an
example of Markov processes enjoying the time-inversion property of
degree 2. The hypothesis of Theorem \ref{thm:main} is indeed
satisfied for $n=\frac{1}{2}m(m+1)$ and
$$
\Phi(x,y) \equiv \(\det(x)\det(y)\)^{-\frac{\delta-m-1}{4}}\
{\bf\tilde I}_{\frac{\delta-m-1}{2}}\(\frac{xy}{4}\),
$$
\begin{equation} \theta
(y) \equiv \frac{1}{2^n} (\det(y))^{\frac{\delta-m-1}{2}}\quad,
\quad\rho(x) \equiv -\frac{1}{2}Tr(x).
\end{equation}

Next we use a skew-product representation, as for the Dunkl process,
to elaborate on the Wishart process and derive a matrix-valued
process with jumps. The skew-product representation allows the
expression of the semigroup densities in terms of the Wishart
transition probability densities.

\begin{definition}
Let $(N_t^{(\lambda)},t\ge0)$ be a Poisson process with intensity
$\lambda$. Let $(X_t,t\ge0)$ be a Wishart process
WIS($\delta,tI_m,\frac{1}{t}x$) independent of the Poisson process.
The skew-Wishart process $(X^{(\lambda)}_t,t\ge0)$ is defined
through the skew-product
\begin{equation}
X^{(\lambda)}_t = X_t\ (-1)^{N_{A_t}^{(\lambda)}}
\end{equation}
where $\displaystyle A_t = \int_0^t Tr(X_s^{-1})ds$.
\end{definition}

\begin{proposition}
The transition probability densities of the skew-Wishart process are
related to the semigroup densities $p_t(x,y)$ of the Wishart process
$X_t$ as follows
\begin{eqnarray}
p^{(\lambda)}_t(x,y) &=& p_t(x,|y|)\left\{1_{\{y\in
S_m^+\}}\frac{1}{2}\(1+\(\frac{{\bf\tilde I_{\nu'}}}{{\bf\tilde
I_\nu}}\)\(\frac{xy}{4t^2}\)\)\right.\nonumber\\
&&\left.\quad\quad\quad+\ 1_{\{y\in
S_m^-\}}\frac{1}{2}\(1-\(\frac{{\bf\tilde I_{\nu'}}}{{\bf\tilde
I_\nu}}\)\(\frac{-xy}{4t^2}\)\)\right\}.
\end{eqnarray}
for $\nu =\frac{\delta-m-1}{2}$, $\nu' = \sqrt{\nu^2+4\lambda}$ and
$|y| = y(1_{\{y\in S_m^+\}}- 1_{\{y\in S_m^-\}})$.
\end{proposition}

\begin{proof}
Let $({\bf P}_t)_{t>0}$ be the semigroup of the skew-Wishart
process. For $x>0$ and $f\in C_c(M_m(\R))$,
\begin{eqnarray*}
{\bf P}_t f(x)& =& {\bf E}_x \[f(X^{(\lambda)}_t)\]\\ &=& {\bf E}_x
\[f(X_t)\ 1_{\{N^{(\lambda)}_{A_t}\ {\rm is\ even}\}}\] + {\bf E}_x
\[f(-X_t)\
1_{\{N^{(\lambda)}_{A_t}\ {\rm is\ odd}\}}\].
\end{eqnarray*}
With $\P(N^{(\lambda)}_{u}\ {\rm is\ even}) =
\frac{1}{2}(1+\exp(-2\lambda u))$, we have
\begin{equation}\label{eq:proofPt}
{\bf P}_t f(x) = {\bf E}_x
\[f(X_t)\ \frac{1}{2}(1+\exp(-2\lambda A_t))\] + {\bf E}_x
\[f(-X_t)\ \frac{1}{2}(1-\exp(-2\lambda A_t))\].
\end{equation}
Let ${\bf Q}_x^{(\nu')}$ with $\nu' = \frac{\delta'-m-1}{2}$ denote
the probability law of a Wishart process WIS$(\delta', t I,
\frac{1}{t} x)$, and ${\bf Q}_x^{(\nu)}$ with $\nu =
\frac{\delta-m-1}{2}$ the probability law of $X_t$. According to
Theorem 1.2 (Remark 2.3) in \cite{DDMY2004}, the probability laws
are related as follows:
$$
{\bf Q}_x^{(\nu')}\big|_{\F_t} = \(\frac{\det X_t}{\det
x}\)^{\frac{\nu'-\nu}{2}} \exp\(-\frac{\nu'^2-\nu^2}{2}\int_0^t
Tr(X_s^{-1})ds\) \cdot {\bf Q}_x^{(\nu)}\big|_{\F_t},
$$
from which we deduce
$$
\frac{{\bf p}_t^{(\nu')}(x,y)}{{\bf p}_t^{(\nu)}(x,y)} =
\(\frac{\det y}{\det x}\)^{\frac{\nu'-\nu}{2}} {\bf
Q}_x^{(\nu)}\[\exp\(-\frac{\nu'^2-\nu^2}{2}\int_0^t
Tr(X_s^{-1})ds\)\bigg| X_t = y\].
$$
Thus, from the expression of the semigroup densities in
(\ref{eq:WishartDensities}), we have
$$
{\bf E}_x^{(\nu)}\[\exp\(-2\lambda\int_0^t Tr(X_s^{-1})ds\)\bigg|
X_t = y\] = \(\frac{{\bf\tilde
I_{\sqrt{\nu^2+4\lambda}}}}{{\bf\tilde I_\nu}}\)\(\frac{xy}{4t^2}\).
$$
Combining the latter with (\ref{eq:proofPt}) yields the semigroup
densities for the skew-Wishart process.
\end{proof}

The skew-Wishart is an example of matrix-valued process with jumps
that enjoys the time-inversion property of degree 2. Indeed, by
setting
\begin{eqnarray*}
\Phi(x,y) &\equiv& \(\det(x)\det(|y|)\)^{-\frac{\nu}{2}}\
\left\{1_{\{y\in S_m^+\}}\frac{1}{2}\({\bf\tilde I}_{\nu}+{\bf\tilde
I_{\nu'}}\)\(\frac{xy}{4}\)\right.\\ && \left.+\ 1_{\{y\in
S_m^-\}}\frac{1}{2}\({\bf\tilde I}_{\nu}-{\bf\tilde
I_{\nu'}}\)\(\frac{-xy}{4}\)\right\},
\end{eqnarray*}
\begin{equation}
\theta (y) \equiv \frac{1}{2^n} (\det(|y|))^{\nu}\quad, \quad\rho(x)
\equiv -\frac{1}{2}Tr(|x|),
\end{equation}
the conditions of Theorem \ref{thm:main} are satisfied for $\alpha =
2$.

\appendix
\section{Generalized hypergeometric functions}
Using the notation in Muirhead \cite{Mu1982}, hypergeometric
functions of matrix arguments are defined for a real symmetric
$m\times m$ matrix $X$, $a_i\in\C$ and $b_j\in\C\backslash
\{0,\frac{1}{2},1,\ldots,\frac{m-1}{2}\}$ by
\begin{equation}
_p{\bf F}_q(a_1, \ldots,a_p;b_1, \ldots,b_q;X) = \sum_{k=0}^\infty
\sum_\kappa \frac{(a_1)_\kappa \cdots(a_p)_\kappa}{(b_1)_\kappa
\cdots(b_q)_\kappa}\frac{C_\kappa(X)}{k!}
\end{equation}
where the second summation is over all partitions $\kappa =
(k_1,\ldots,k_m),\ k_1\ge\cdots\ge k_m\ge0$, of $k=\sum_{i=1}^m
k_i,\ k! = k_1!\cdots k_m!$ and the generalized Pochhammer symbols
are given by
$$
(a)_\kappa= \prod_{i=1}^m \(a-\frac{i-1}{2}\)_{k_i},\quad (a)_k =
a(a+1)\cdots(a+k-1), \quad (a)_0=1.
$$
$C_\kappa(X)$ is the zonal polynomial corresponding to $\kappa$,
which is a symmetric, homogeneous polynomial of degree $k$ in the
eigenvalues of $X$ that satisfies
\begin{equation}
C_\kappa(YX) = C_\kappa(\sqrt YX\sqrt Y)
\end{equation}
for some $Y\in S_m^+$. The function $_p{\bf F}_q(a_1,
\ldots,a_p;b_1, \ldots,b_q;YX)$ thus makes sense. Finally, we define
the generalized modified Bessel function by
\begin{equation}
{\bf\tilde I_\nu}(X) =
\frac{\(\det(X)\)^{\frac{\nu}{2}}}{{\bf\Gamma}_m\(\nu
+\frac{m+1}{2}\)}\ {}_0{\bf F}_1\(\nu+\frac{m+1}{2};X\)
\end{equation}
where the generalized gamma function is given as a product of the
usual gamma functions,
\begin{equation}
{\bf\Gamma}_m\(\alpha\)= \pi^{\frac{m(m-1)}{4}} \prod_{i=1}^m
\Gamma\(\alpha-\frac{i-1}{2}\)
\end{equation}
for Re$(\alpha)>\frac{m-1}{2}$. Note that the generalized modified
Bessel for $m=1$ relates to the usual modified Bessel function
$I_\nu(x)$ by ${\bf\tilde I_\nu}(x) = I_\nu(2\sqrt x)$ (see
\cite{Le1972}).

\section*{Acknowledgements} The author wishes to thank Marc Yor for
helpful discussions and references.


\begin{thebibliography}{99}
\bibitem{Br1989a}
      Bru, M.F.,
     \emph{Processus de {W}ishart},
      C.R. Acad. Sci. Paris, S{\'e}rie I,
     {\bf t.308},
     pp. 29--32,
      1989

\bibitem{Br1989b}
      Bru, M.F.,
     \emph{Diffusions of perturbed principal component analysis},
      J. Multivariate Anal.,
     {\bf 29},
     pp. 127--136,
      1989

\bibitem{Br1991}
      Bru, M.F.,
     \emph{Wishart processes},
      J. Theo. Probab.,
     {\bf 4},
     pp. 725--751,
      1991

\bibitem{CL2001}
      C{\'e}pa, E. and L{\'e}pingle, D.,
     \emph{Brownian particles with electrostatic repulsion on the circle:
                    {D}yson's model for unitary random matrices revisited},
      ESAIM Probability and Statistics,
     {\bf 5},
     pp. 203--224,
      2001

\bibitem{CGM2005}
      Chouchene, F. and Gallardo, L. and Mili, M.,
     \emph{The heat semigroup for the {J}acobi-{D}unkl
                    operator and the related {M}arkov processes},
      Pot. Anal.,
      {\bf 25(2)},
      pp. 103--119,
      2006

\bibitem{Ch2005}
      Chybiryakov, O.,
     \emph{Skew-product representations of multidimensional
                    {D}unkl {M}arkov processes},
      Preprint,
      2005

\bibitem{DDMY2004}
      Donati-Martin, C. and Doumerc, Y. and Matsumoto, H. and Yor, M.,
     \emph{Some properties of the {W}ishart processes and a matrix
                extension of the {H}artman-{W}atson laws},
      Publ. RIMS, Kyoto Univ.,
     {\bf 40},
     pp. 1385--1412,
      2004

\bibitem{Du1989}
      Dunkl, C.,
     \emph{Differential-difference operators associated
                    to reflection groups},
      Trans. Amer. Math. Soc.,
     {\bf 311(1)},
     pp. 167--183,
      1989

\bibitem{Du1992}
      Dunkl, C.,
     \emph{Hankel transforms associated to finite reflection
                    groups},
      Contemp. Math.,
     {\bf 138},
     pp. 123--138,
      1992

\bibitem{Dy1975}
      Dynkin, E.B.,
     \emph{Markov representations of stochastic systems},
      Russ. Math. Surv.,
     {\bf 30}(1),
     pp. 65--104,
      1975

\bibitem{Dy1962}
      Dyson, F.J.,
     \emph{A {B}rownian-motion model for the eigenvalues of a
                    random matrix},
      J. Math. Phys.,
     {\bf 3},
     pp. 1191--1198,
      1962

\bibitem{GY2005a}
      Gallardo, L. and Yor, M.,
     \emph{Some remarkable properties of the {D}unkl martingales},
      S\'em. Proba. XXXIX, Lect. Notes in Math,
      {\bf 1874},
      2006

\bibitem{GY2005b}
      Gallardo, L. and Yor, M.,
     \emph{Some new examples of {M}arkov processes which enjoy
                    the time-inversion property},
      Probab. Theory Relat. Fields,
     {\bf 132},
     pp. 150--162,
      2005

\bibitem{GY2004c}
      Gallardo, L. and Yor, M.,
     \emph{A chaotic representation property of the
                    multidimensional {D}unkl processes},
      Ann. Prob.,
      {\bf 34(4)},
      pp. 1530--1549,
      2006

\bibitem{Gr1999}
      Grabiner, D.,
     \emph{Brownian motion in a {W}eyl chamber, non-colliding
                    particles, and random matrices},
      Annales de l'I. H. P. Probabilit\'es et Statistiques,
     {\bf 35},
     pp. 177--204,
      1999

\bibitem{G1997}
      Gruet, J.C.,
     \emph{Windings of hyperbolic {B}rownian motion,
                    in {E}xponential functionals and principal
                    values related to {B}rownian motion},
      Bibl. Rev. Math. Ibroamericana,
                    Rev. Math. Ibroamericana, Madrid,
     pp. 35--72,
      1997

\bibitem{G1998}
      Gruet, J.C.,
     \emph{Jacobi radial stable processes},
      Ann. Math. Blaise Pascal,
     {\bf 5},
     pp. 39--48,
      1998

\bibitem{IM1974}
      Ito, K. and McKean, H.P.,
     \emph{Diffusion Processes and Their Sample Paths},
      Springer-Verlag, Berlin,
      1974

\bibitem{KO2001}
      K{\"o}nig, W. and O'Connell, N.,
     \emph{Eigenvalues of the {L}aguerre process as non-colliding
                    squared {B}essel processes},
      Elect. Comm. in Probab.,
     {\bf 6},
     pp. 107--114,
      2001

\bibitem{Ku1977}
      Kuznetsov, S.E.,
     \emph{Construction of a regular split process},
      Theor. Prob. Appl.,
     {\bf 22},
     pp. 773--793,
      1977

\bibitem{Mu1982}
      Muirhead, R.J.,
     \emph{Aspects of Multivariate Statistic Theory},
      Wiley, New York,
      1982

\bibitem{La1972}
      Lamperti, J.,
     \emph{Semi-Stable {M}arkov Processes {I}},
      Zeit. f{\"u}r Wahr.,
     {\bf 22},
     pp. 205--255,
      1972

\bibitem{Le1972}
      Lebedev, N.N.,
     \emph{Special Functions and Their Applications},
      Dover, New York,
      1972

\bibitem{PY1981}
      Pitman, J. and Yor, M.,
     \emph{Bessel processes and infinitely divisible laws},
      In: Stochastic Integrals. ed. D. Williams.,
                    Lect. Notes in Math,
     {\bf 851},
      1981

\bibitem{Ro1998}
      R{\"o}sler, M.,
     \emph{Generalized {H}ermite polynomials and the heat
                    equations for {D}unkl operators},
      Comm. Math. Phys.,
     {\bf 192},
     pp. 519--542,
      1998

\bibitem{RV1998}
      R{\"o}sler, M. and Voit, M.,
     \emph{Markov processes related with {D}unkl operators},
      Adv. in App. Math.,
     {\bf 21(4)},
     pp. 575--643,
      1998

\bibitem{Wa1975}
      Watanabe, S.,
     \emph{On time inversion of one-dimensional diffusion
                    processes},
      Z. Wahrsch. Verw. Gebiete,
     {\bf 31},
     pp. 115--124,
      1975

\end{thebibliography}
\end{document}